\documentclass[10pt]{amsart}
\usepackage[utf8]{inputenc}
\usepackage[english]{babel}
\usepackage{amsmath}
\usepackage{amssymb}
\usepackage{amsthm}
\usepackage{latexsym}
\usepackage{amsfonts}
\usepackage{bbm}
\usepackage{mathrsfs}
\usepackage{textcomp}
\usepackage[T1]{fontenc}
\usepackage{graphicx}
\usepackage{setspace}
\usepackage{nicefrac}
\usepackage{indentfirst}
\usepackage{enumerate}
\usepackage{tikz}
\usepackage{wasysym}
\usepackage{upgreek}
\usepackage{hyperref}
\usepackage{booktabs}
\usepackage{color}

\usepackage[labelfont=bf, font={small}]{caption}[2005/07/16]

\usepackage{epigraph}
\usepackage[a4paper,top=3cm,bottom=3cm,left=3cm,right=3cm,bindingoffset=1cm]{geometry}

\newtheorem{teo}{Theorem}[section]
\newtheorem{prop}[teo]{Proposition}
\newtheorem{lem}[teo]{Lemma}
\newtheorem{corol}[teo]{Corollary}
\newtheorem{congettura}[teo]{Conjecture}

\theoremstyle{definition}
\newtheorem{defin}[teo]{Definition}

\theoremstyle{remark}
\newtheorem{remark}[teo]{Remark}
\newtheorem{notation}[teo]{Notation}
\newtheorem*{PROOF1.1}{Proof of Theorem 1.1}

\numberwithin{equation}{section}

\title{An asymptotic bound for secant varieties of Segre varieties}
\author{Fulvio Gesmundo}

\newcommand{\PP}{\mathbbm{P}}

\newcommand{\RRR}{\mathscr{R}}
\newcommand{\OOO}{\mathscr{O}}

\newcommand{\ootimes}{\otimes \dots \otimes}
\newcommand{\vvirg}{ , \dots , }
\newcommand{\ttimes}{ \times \dots \times }

\begin{document}


\begin{abstract}
This paper studies the defectivity of secant varieties of Segre varieties. We prove that there exists an asymptotic lower estimate for the greater non-defective secant variety (without filling the ambient space) of any given Segre variety. In particular, we prove that the ratio between the greater non-defective secant variety of a Segre variety and its expected rank is lower bounded by a value depending just on the number of factors of the Segre variety. Moreover, in the final section, we present some results obtained by explicit computation, proving the non-defectivity of all the secant varieties of Segre varieties of the shape $\left(\PP^{n}\right)^4$, with $2 \leq n\leq 10$, except at most $\sigma_{199}(\left(\PP^8\right)^4)$ and $\sigma_{357}(\left(\PP^{10}\right)^4)$.
\end{abstract}

\maketitle

\section{Introduction}

Let $X \subseteq \PP V$ be a projective variety of dimension $d$, where $\PP V$ denotes the projective space of lines of an $(n+1)$-dimensional vector space $V$ over an algebraically closed field of characteristic zero $K$. $\sigma_s(X)$ denotes the $s$-secant variety of $X$, namely the Zariski closure of the union of the linear span of $s$-tuples of general points $p_1 \vvirg p_s \in X$. More precisely

$$\sigma_s(X) = \overline{\bigcup_{\scriptsize \begin{array}{c}p_i \in X\\ i=1 \vvirg s \end{array}} \PP^{s-1}_{p_1 \dots p_s}},$$

where $\PP^{s-1}_{p_1 \dots p_s}$ denotes the $(s-1)$-dimensional projective space generated by the points $p_1 \vvirg p_s$.

The maximum possible dimension of $\sigma_s(X)$ is $\min \left\{ n, s(d+1) -1 \right\}$ and this value is said the expected dimension of the secant variety. A secant variety is said defective if its dimension is lower than the expected dimension; otherwise it is said non-defective.

For $i=1 \vvirg k$, let $V_i$ be a $(n_i+1)$-dimensional vector space over $K$. We denote by $\PP^{n_i}$ the projective space $\PP V_i$. The Segre variety of the spaces $V_1 \vvirg V_k$ is defined to be the variety $X=\PP^{n_1} \ttimes \PP^{n_k}$, embedded in $\PP(V_1 \ootimes V_k)$ by the Segre embedding. The expected dimension of $\sigma_s (X)$ is 
\begin{equation}\label{eq:expecteddimension}
\min \left\{\prod (n_i+1) -1 , s(1 + \sum n_i) -1 \right\}.
\end{equation}

The interest of this subject is due to the link between this topic and the study of tensor rank. The image of the Segre embedding (i.e. the embedded Segre variety) is viewed as the set of rank one tensors, or decomposable tensors. A tensor $\omega$ is said to have \emph{rank} $r$, if $r$ is the smallest positive integer such that $\omega$ can be written as the sum of $r$ decomposable tensors. A tensor $\omega$ is said to have \emph{border rank} $r$, if $r$ is the smallest positive integer such that $\omega$ is the limit of a sequence of rank $r$ tensors. With this notation, the $s$-secant variety is viewed as the set of tensors whose border rank is at most $s$.

Thus, our aim is to study the defectivity of secant varieties of a Segre variety $X$, in order to determine the minimum integer $R$ such that all the tensors in $V_1 \ootimes V_k$ have at most border rank $R$; equivalently, $R$ is the smallest positive integer such that $\sigma_R (X)$ fills the ambient space. In accordance with \cite{burg}, \cite{comon} et al., this value is said the generic rank of $X$ and is denoted by $\underline{\mathbf{R}}(X)$. 

For example, if a given Segre variety $\PP V_1 \ttimes \PP V_k$ has no defective secant varieties, i.e. all the secant varieties have the expected dimension, then we can deduce that the generic rank is the smallest integer $s$ such that the minimum in \eqref{eq:expecteddimension} is given by the left term. This value is 
$$R = \left\lceil \frac{\prod(n_i+1)}{1+\sum n_i} \right\rceil.$$
Thus, the $R$-secant variety is the first one which fills the ambient space, and all the tensors in $V_1 \ootimes V_k$ have border rank $\leq R$.

In Section 2, we recall some preliminary results about this topic. In particular, in \cite{aboottapet} we find an inductive technique that makes it possible to reduce the study of Segre varieties of high dimension to the study of Segre varieties of lower dimension. We will sum up this procedure by the language of the reduction trees.

In Section 3, we introduce the notions of room and safety region: in this section we give a lower bound for the positive safety bound of any Segre variety and an upper bound for Segre varieties of the shape $\left(\PP^1\right)^k$; moreover we provide a result that will help us to compute reduction trees.

In Section 4, we provide some results that make it possible to control the growth of the safety bounds as the dimensions of the factors of the Segre variety change. Given a Segre variety $X$ with $k$ factors, these results allow us to prove that the ratio between the greatest $s$ such that the $s$-secant variety is non-defective without filling the ambient space and the expected border rank of $X$ is a constant $\Theta_k$ depending just on $k$:
\begin{teo}\label{Thm: stima asintotica per la non difettivita}
 Let $k\geq 3$ be a positive integer and let $X=\PP^{n_1} \ttimes \PP^{n_k}$ be a Segre variety with $k$ factors. Then there exists a coefficient $\Theta_k \in (0,1)$, depending just on $k$, such that the $s$-secant variety of $X$ is non-defective for any $s$ such that
$$ s \leq \Theta _k \frac{\prod_1 ^k (n_i +1)}{1+\sum_1 ^k n_i}.$$
\end{teo}
The coefficient $\Theta_k$ is the ``ratio'' between the number of values for $s$ where we can prove that $\sigma_s(X)$ is non defective and the number of possible values for $s$. In Corollary \ref{cor: Teo 1.1} we prove that in case $n_i+1$ are powers of $2$ then $\Theta_k\to 1$ when $k\to +\infty$.

Section 5 consists in the proof of this theorem.

Section 6 has to do with results obtained by the aid of computer: we compute the safety region of low-dimensional Segre varieties, remarking an unexpected symmetry between the positive bound and the negative one. Then, we sum up all the information collected about the safety bounds of a Segre variety. Moreover, dealing with low-dimensional varieties of the shape $\left( \PP^n \right)^4$, we prove that these are non-defective if $2\leq n\leq7$ or $n=9$ and have at most one defective secant variety if $n=8$ or $n=10$; this supports the suspect claimed in \cite{aboottapet}, according to which the only defective secant varieties of Segre varieties with this shape are $\PP^2\times \PP^2 \times \PP^2$ and $\PP^1 \times \PP^1 \times \PP^1 \times \PP^1$.

This paper is based on my Master's thesis and on the results I achieved in the following few months. I would like to express my deepest gratitude to Giorgio Ottaviani for his invaluable advice and for all the efforts he spent for me in this period.

\section{Preliminaries}

By an easy application of the Leibniz product rule (see \cite{aboottapet}, Lemma 2.2), we deduce that the affine tangent space to a Segre variety $X=\PP^{n_1} \ttimes \PP^{n_k}$ at a point $p=v_1 \ootimes v_k$ is
\begin{equation*}
 T_p X = V_1 \otimes v_2 \ootimes v_k + \dots + v_1 \ootimes v_{k-1} \otimes V_k.
\end{equation*}
In order to lighten the notation, often we will omit the symbol of tensor product, by denoting, for example, a point $u \otimes v$ simply by $uv$. Moreover, we use the notation $(T_p X)_j$, referring to the $j$-th summand of the tangent space to the Segre variety at $p = v_1 \dots v_k$, namely:
$$(T_p X)_j = v_1 \dots v_{j-1}V_j v_{j+1} \dots v_k .$$

We recall the Terracini's Lemma, that makes it possible to compute the tangent space of a secant variety of a variety $X$ at a general point, by computing the linear span of the tangent spaces of $X$. More precisely:
\begin{teo}[Terracini's Lemma]
 Let $X \subseteq \PP^n$ be a projective variety and let $p_i$ be general points of $X$, for $i=1 \vvirg s$. Then, set $q$ be a general point of $\sigma_s(X)$ such that $q \in \PP^{s-1}_{p_1 \vvirg p_s}$, we have:

$$T_q \sigma_s(X) = T_{p_1} X + \dots + T_{p_s} X.$$
\end{teo}

Moreover, we recall the main result of \cite{aboottapet}, providing an inductive technique that reduces the computation of $\dim \sigma_s (X)$ to the computation of the dimension of partial secant varieties to lower dimensional Segre varieties. We introduce some definitions and notations.
\\

\begin{notation}
 Let $X = \PP ^{n_1} \ttimes \PP^{n_k}$. Sometimes, we can denote this Segre variety by $\PP^{\mathbf{n}}$, where $\mathbf{n}=(n_1 \vvirg n_k)$. For $s$ general points $p_1 \vvirg p_s \in X$, $T_s X$ denotes the space $T_{p_1} X + \dots + T_{p_s}X$. Analogously, $(T_s X)_j$ denotes the space $(T_{p_1} X)_j + \dots +(T_{p_s} X)_j$.
\end{notation}

\begin{defin}\label{def:Tnsa}
Let $X = \PP^{n_1} \ttimes \PP^{n_k}$ be a Segre variety and let $s,a_1 \vvirg a_k$ be non-negative integers. We will say that $T(n_1 \vvirg n_k;s;a_1 \vvirg a_k)$ (sometimes simply $T(\mathbf{n},s,\mathbf{a})$) is true if the space
$$L = T_s X + (T_{a_1} X)_1 + \dots +(T_{a_k} X)_k$$
has dimension exactly $D$, where
$$D= \min \left\{s\left( 1+\textstyle \sum _{1} ^k n_i \right) + \textstyle\sum_1 ^k a_i(n_i +1) , \textstyle\prod_1 ^k (n_i +1) \right\}.$$

\begin{itemize}
 \item[-] $T(\mathbf{n},s,\mathbf{a})$ is said \emph{subabundant} if the minimum is verified by the first term, namely if $s\left( 1+\textstyle \sum _{1} ^k n_i \right) + \sum_1 ^k a_i(n_i +1) \leq \prod_1 ^k (n_i +1)$. 
\item[-] $T(\mathbf{n},s,\mathbf{a})$ is said \emph{superabundant} if the inequality holds in the opposite direction, namely if $s\left( 1+\textstyle \sum _{1} ^k n_i \right) + \sum_1 ^k a_i(n_i +1) \geq \prod_1 ^k (n_i +1)$.
\item[-] $T(\mathbf{n},s,\mathbf{a})$ is said \emph{equiabundant} if equality holds, namely if it is both subabundant and superabundant.
\end{itemize}
\end{defin}

We recall some useful remarks, proved in \cite{aboottapet}.

\begin{remark}\label{remark: ordinestatements}
By Terracini's Lemma, we can deduce that the space $L$ of Definition \ref{def:Tnsa} is the affine tangent space to a partial secant variety to $X$ and that $D$ is its expected dimension (see \cite{aboottapet}, Remark 3.10). Thus, we may apply semicontinuity arguments, using specialized points in order to prove that expected dimension holds in the general case.

Moreover, given non-negative integers $n_1 \vvirg n_k,n_1 ' \vvirg n_k ' , s,s' ,a_1 \vvirg a_k,a_1' \vvirg a_k'$ such that $n_i \leq n_i '$,$s\leq s',a_j \leq a_j'$ for all $i$ and $j$, we make the following remarks (\cite{aboottapet}, Remark 3.3 and Theorem 3.11):
\begin{itemize}
\item[-] If $T(\mathbf{n},s',\mathbf{a}')$ is true and subabundant, then $T(\mathbf{n},s,\mathbf{a})$ is true and subabundant.
\item[-] If $T(\mathbf{n},s,\mathbf{a})$ is true and superabundant, then $T(\mathbf{n},s',\mathbf{a}')$ is true and superabundant.
\item[-] If $T(\mathbf{n},s,\mathbf{a})$ is true and subabundant, then $T(\mathbf{n}',s,\mathbf{a})$ is true and subabundant.
\item[-] If $T(\mathbf{n}',s,\mathbf{a})$ is true and superabundant, then $T(\mathbf{n},s,\mathbf{a})$ is true and superabundant.
\end{itemize}
\end{remark}

Now we can state the theorem that provides the inductive technique introduced in \cite{aboottapet}.

\begin{teo}\label{teosub teosuper}
 Let $n_k +1 = (n_k '+1) + (n_k ''+1), s=s'+s'', a_j = a_j' + a_j''$ for $j=1 \vvirg k-1$. Suppose
\begin{enumerate}[(i)]
 \item $T(n_1 \vvirg n_{k-1}, n_k';s';a_1 ' \vvirg a_{k-1}', a_k +s'')$ is true and subabundant (resp. superabundant),
\item $T(n_1 \vvirg n_{k-1}, n_k'';s'';a_1 '' \vvirg a_{k-1}'', a_k +s')$ is true and subabundant (resp. superabundant).
\end{enumerate}
Then $T(n_1 \vvirg n_k;s;a_1 \vvirg a_k)$ is true and subabundant (resp. superabundant).
\end{teo}

Thus, in order to study whether a given $T(\mathbf{n},s,\mathbf{a})$ is true, we can reduce this statement into two statements in accordance to the arithmetic relations of Theorem \ref{teosub teosuper}, and verify whether these statements are true. The two statements in the reduction describe the defectivity of two partial secant varieties of two Segre varieties of dimension lower than $\PP^{\mathbf{n}}$. Thus, our strategy consists in iterating the application of Theorem \ref{teosub teosuper}, in order to obtain statements that can be checked by explicit computation, as described in \cite{aboottapet}.

Furthermore, we recall the following useful result, that allows us to cut out a factor when its dimension is zero (see \cite{aboottapet}, Corollaries 3.8 and 3.9):
\begin{prop}\label{prop:fattorizerotagliati}
 Let $n_2 \vvirg n_k,s,a_1 \vvirg a_k$ be nonnegative integers. Then, the following statements hold:
\begin{enumerate}
\item if $T(n_2 \vvirg n_k; s; a_2 \vvirg a_k)$ is true, then $T(0,n_2 \vvirg n_k; s; a_1,a_2 \vvirg a_k)$ is true;
\item if $T(0,n_2 \vvirg n_k;s;a_1,a_2 \vvirg a_k)$ is true and subabundant, then $T(n_2 \vvirg n_k ; s ;a_2 \vvirg a_k)$ is true and subabundant;
\item if $T(0,n_2 \vvirg n_k;s;0,a_2 \vvirg a_k)$ is true and superabundant, then $T(n_2 \vvirg n_k ; s ;a_2 \vvirg a_k)$ is true and superabundant.
\end{enumerate}
\end{prop}

\begin{defin}
 Let $\mathcal{T} = T(\mathbf{n},s,\mathbf{a}), \mathcal{T} = T(\mathbf{n}',s',\mathbf{a}'), \mathcal{T} = T(\mathbf{n}'',s'',\mathbf{a}'')$ be three statements. We will say that $\mathcal{T}$ is reduced to $\mathcal{T}'$ and $\mathcal{T}''$ on the $j$-th factor if these two statements satisfy the arithmetic relations of Theorem \ref{teosub teosuper} with respect of the $j$-th factor, namely $n_i ' = n_i '' = n_i$ for $i \neq j$, $(n_j'+1)+(n_j'' +1)= n_j +1$, $s'+s''=s$, $a_i' +a_i'' = a_i$ for $i\neq j$, $a_j' = a_j+s''$ and $a_j'' = a_j +s'$. Moreover, if abundance is preserved, i.e. if the statements are all subabundant or all superabundant, we say that the reduction is \emph{eligible}.
\end{defin}

Often it will be useful to imagine one or more reductions in a binary tree, where a node $\mathcal{T}$ has two children $\mathcal{T}'$ and $\mathcal{T}''$ if the statement $\mathcal{T}$ is reduced in $\mathcal{T}'$ and $\mathcal{T}''$. A tree that describes the reduction for a statement $T(\mathbf{n},s,\mathbf{a})$ is said to be a \emph{reduction tree} for $T(\mathbf{n},s,\mathbf{a})$.

\section{Room and safety region}

In order to study how we can determine eligible reduction trees for a generic statement $T(\mathbf{n},s,\mathbf{a})$, we introduce a fundamental notion, which gives a sort of measure of ``how much'' a statement is subabundant or superabundant.

\begin{defin}
 Consider a statement $T(\mathbf{n},s,\mathbf{a})$. We call \emph{room} of $T(\mathbf{n},s,\mathbf{a})$ the value
$$ \RRR = \prod_{i=1}^k (n_i +1) - \left[ s\left( 1 + \sum_{i=1}^k n_i\right) + \sum_{i=1}^k a_i(n_i+1)\right].$$
\end{defin}

When we need to focus on a specific factor of the Segre variety, for example on the $j$-th one, we will use the following equivalent form of the room:

\begin{equation}\label{eq: room j-evidenza}
\RRR = (n_j +1) \left[ \prod_{i\neq j} (n_i+1) - (s +a_j) \right] -s \left[ \sum_{i\neq j} n_i\right] - \sum_{i \neq j} a_i (n_i+1). 
\end{equation}

Obviously, we have that $T(\mathbf{n},s,\mathbf{a})$ is subabundant if and only if its room is nonnegative and $T(\mathbf{n},s,\mathbf{a})$ is superabundant if and only if its room is non positive.

Moreover, it is easy to check that the value of room is additive in a reduction, namely we have:

\begin{prop}\label{prop: room additiva}
 Let $T(\mathbf{n},s,\mathbf{a})$ be a statement with room $\RRR$. Let $T(\mathbf{n}',s',\mathbf{a}')$ and $T(\mathbf{n}'',s'',\mathbf{a}'')$ be two statements in which $T(\mathbf{n},s,\mathbf{a})$ can be reduced, and let $\RRR'$ and $\RRR''$ be their rooms. Then, $\RRR = \RRR' + \RRR''$.
\begin{proof}
 Suppose that the reduction is applied on the $j$-th factor. Consider the rooms of the statement, written as in \eqref{eq: room j-evidenza}:
\begin{equation*}
\begin{split}
\RRR &= (n_j +1) \left[\textstyle \prod_{i\neq j} (n_i+1) - (s +a_j) \right] -s \left[ \textstyle\sum_{i\neq j} n_i\right] -\textstyle \sum_{i \neq j} a_i (n_i+1),  \\
\RRR' &= (n_j' +1) \left[ \textstyle \prod_{i\neq j} (n_i+1) - (s' +a_j') \right] -s' \left[ \textstyle\sum_{i\neq j} n_i\right] - \textstyle \sum_{i \neq j} a_i' (n_i+1), \\
\RRR'' &= (n_j'' +1) \left[ \textstyle\prod_{i\neq j} (n_i+1) - (s'' +a_j'') \right] -s'' \left[\textstyle \sum_{i\neq j} n_i\right] -\textstyle \sum_{i \neq j} a_i'' (n_i+1).
\end{split}
\end{equation*}
Since $s' +a_j' = s'' + a_j'' = s + a_j$, and $n_j+1 = (n_j'+1)+(n_j''+1)$, $s=s'+s''$ and $a_i=a_i'+a_i''$ for $i\neq j$, we conclude.
\end{proof}
\end{prop}

Our aim is to provide some results which permit us to control how the room is distributed in a reduction. In other words, we search for conditions that allow us to maximize (in absolute value) the value of room of $\mathcal{T}'$ and $\mathcal{T}''$, while reducing the statement $\mathcal{T}$. Thus, consider a statement $T(n_1 \vvirg n_k;s;a_1 \vvirg a_k)$ to be reduced on the $k$-th factor. The integer $n_k +1$ is divided into two integers $n_k'+1$ and $n_k''+1$. We can consider the ratio between $n_k' +1$ and $n_k +1$, which is a value lower than $1$. By examples and by tests that we can run by the aid of a computer, it seems that a good strategy to find optimal reductions is to look for a reduction in which the room is divided proportionally to the division of $n_k +1$. The following result shows that, if some conditions hold, we can find a reduction such that the difference between its room and the optimal one is limited by a term depending just on the dimensions of the factors of the Segre variety.

\begin{teo}\label{teo: room proporzionata}
 Let $k \geq 3$, let $n_1 \vvirg n_k,s,a_1 \vvirg a_k$ be positive integers and let $\mathcal{T}=T(n_1 \vvirg n_k;s;a_1 \vvirg a_k)$. Let $\RRR$ be the room of $\mathcal{T}$. Let $n_k ' , n_k''$ be positive integers such that $n_k+1= (n_k' +1)+(n_k'' +1)$ and let $\gamma' = \frac{n_k'+1}{n_k+1}$ and $\gamma''= \frac{n_k''+1}{n_k+1}$. Then, there exists a reduction of $\mathcal{T}$ on the $k$-th factor into two statements $\mathcal{T}'$ and $\mathcal{T}''$, whose rooms are $\RRR'$ and $\RRR''$ respectively, such that:

\begin{equation}\label{eq:condizioneroomproporzionata}
\begin{split}
\left\vert \RRR' - \gamma' \RRR \right\vert \leq \frac{1}{2} \sum_{i=1}^{k-1} n_i , \\
\left\vert \RRR'' - \gamma'' \RRR \right\vert \leq \frac{1}{2} \sum_{i=1}^{k-1} n_i .
\end{split}
\end{equation}

Moreover, if we have

\begin{equation}\label{eq:condizioneroomproporzionataammissibile}
\frac{1}{2} \sum _{i=1} ^{k-1} n_i \leq \min \{ \gamma', \gamma'' \} \vert \RRR \vert , 
\end{equation}

then this reduction is eligible.
\begin{proof}
Note that the two inequalities in \eqref{eq:condizioneroomproporzionata} are equivalent by Proposition \ref{prop: room additiva}, then we will prove just the first one. Moreover, if $s$ and all the $a_i$ except $a_k$ are zero, we can easily conclude, since the reduction is determined just by the choice of $n_k'$, and satisfies $\RRR'=\gamma' \RRR$ and $\RRR'' = \gamma'' \RRR$. Thus, suppose that at least one in the set $\{s,a_1 \vvirg a_{k-1} \}$ is other than zero.

First, consider the case $s\geq 1$. Fix the values $a_1 ' \vvirg a_{k-1}'$ (and then $a_1'' \vvirg a_{k-1}''$) and let $\underline{a}'$ denote the vector $(a_1 ' \vvirg a_{k-1}')$. For any $s'=0 \vvirg s$, let $\RRR_{s'}^{\underline{a}'}$ denote the room of the statement $T(n_1 \vvirg n_{k-1},n_k';s';\underline{a}',a_k+s-s')$. We have
\begin{equation} \label{eq:teo room proporzionata - distanza termini}
\RRR^{\underline{a}'} _{s'} - \RRR^{\underline{a}'} _{s'+1} = \sum_{i=1}^{k-1} n_i . 
\end{equation}
We prove that there is a vector $\underline{a}'$ such that 
\begin{equation}\label{eq:teo room proporzionata - attraversare gammaR}
\RRR ^{\underline{a}'} _0 \geq \gamma ' \RRR \geq \RRR ^{\underline{a}'} _s .
\end{equation}

Arrange all the suitable $\underline{a}'$ in three sets:
\begin{itemize}
\item[-] $\underline{a}'$ such that $\RRR^{\underline{a}'} _{s} > \gamma' \RRR$, and then $\RRR^{\underline{a}'} _{s'} > \gamma' \RRR$ for any $s'$;
\item[-] $\underline{a}'$ such that $\RRR^{\underline{a}'} _{0} < \gamma' \RRR$, and then $\RRR^{\underline{a}'} _{s'} < \gamma' \RRR$ for any $s'$;
\item[-] $\underline{a}'$ such that \eqref{eq:teo room proporzionata - attraversare gammaR} holds.
\end{itemize}
Supposing that the third set is empty, note that the other two sets are non-empty: since $\RRR^{\underline{0}} _{0} > \gamma'\RRR$ and $\RRR^{\underline{a}} _{s} < \gamma' \RRR$, we have that $\underline{0}$ is in the first set and $\underline{a}$ is in the second one.

Now, consider $\tilde{\underline{a}}'$ in the second set which is minimal, namely such that for any vector $\epsilon$ such that $\tilde{\underline{a}}' - \epsilon \geq \underline{0}$, we have that $\tilde{\underline{a}}' - \epsilon$ is in the first set.

Let $\bar{\jmath}$ be one of the indexes such that $a_{\bar{\jmath}} \neq 0$ and let $\tilde{\epsilon}=(0 \vvirg 0,1,0\vvirg 0)$, where $1$ is at the $\bar{\jmath}$-th entry. By the minimality of $\tilde{\underline{a}}$ we have

\begin{equation*}
\begin{split}
\RRR ^{\underline{\tilde{a}}'} _0 < \gamma ' \RRR, \\
\RRR ^{\underline{\tilde{a}}' - \tilde{\epsilon}} _s > \gamma ' \RRR.
\end{split}
\end{equation*}

By subtraction we obtain

$$ s \sum_{i=1}^{k-1} n_i - (n_{\bar{\jmath}} +1) < 0,$$

that is a contradiction since $s \geq 1$ and $k \geq 3$. Thus, the set of $\underline{a}'$ such that \eqref{eq:teo room proporzionata - attraversare gammaR} holds is non-empty. Considering such $\underline{a}'$, we can find $s'$ such that 
\begin{equation*}
\left\vert \RRR ^{\underline{a}'} _{s'} - \gamma' \RRR \right\vert \leq \frac{1}{2} \sum_{i=1} ^{k-1} n_i.
\end{equation*}

In order to settle the case $s=0$, let $d$ be an index such that $a_d \neq 0$ and $n_d$ is maximum among the integers $n_j$ such that $a_j \neq 0$. By the same argument used in the previous part of the proof, we can consider $\underline{a}'= (a_1' \vvirg a_{d-1},a_{d+1} \vvirg a_{k-1})$ and, letting $a_d'$ vary from $0$ to $a_d$ (as in the previous part $s'$ varied from $0$ to $s$), we can show that there exists $\underline{a}'$  such that
$$ \left\vert \RRR' - \gamma ' \RRR \right\vert \leq \frac{1}{2} (n_{\bar{\jmath}} +1), $$
where $n_{\bar \jmath}$ is maximum among all the $a_j$ with $j \neq k$.

This concludes the main part of the proof.

Finally, condition \eqref{eq:condizioneroomproporzionataammissibile} guarantees that the obtained statements have the same abundance, i.e. that the reduction is eligible.
\end{proof}
\end{teo}

We will apply Theorem \ref{teo: room proporzionata} mainly in its first part, without initially caring about the eligibility of the reduction. Indeed, condition \eqref{eq:condizioneroomproporzionataammissibile} is quite limiting and awkward to be used. Thus, we prefer to apply Theorem \ref{teo: room proporzionata} just to find a reduction suitable for our intents, and each time to find an `ad hoc' condition which guarantees its eligibility.

The following remark allows us to define the safety region of a Segre variety:

\begin{remark}
 Fix a Segre variety $\PP^{\mathbf{n}}$. The number of false statements $T(\mathbf{n},s,\mathbf{a})$ is finite. Obviously all the subabundant statements are finite, then the false ones are finite too. As for the superabundant statements, consider that, since the points used for defining the space $L$ in Definition \ref{def:Tnsa} are in general position, the dimension of this space is at least $s + \sum a_i$.
\end{remark}

\begin{defin}
 Let $X = \PP^{n_1} \ttimes \PP^{n_k}$ be a Segre variety. Define
\begin{equation*}
\begin{split}
 \RRR_+ &= \max \left\{ \RRR \text{ such that $\RRR$ is room of a false subabundant statement $T(\mathbf{n},s,\mathbf{a})$} \right\}, \\
 \RRR_- &= \min \left\{ \RRR \text{ such that $\RRR$ is room of a false superabundant statement $T(\mathbf{n},s,\mathbf{a})$} \right\}.
\end{split}
\end{equation*}
Set $\OOO _+ = \RRR_+ +1$ and $\OOO_- = \RRR_- -1$. We call $\OOO_+$ and $\OOO_-$ the \emph{safety bounds} of $X$. The set

$$ (- \infty , \OOO_-] \cup [\OOO_+ ,+\infty) $$

is said the \emph{safety region} of $X$.
\end{defin}

Note, as trivial consequence of the definition, that a statement $T(\mathbf{n},s,\mathbf{a})$ whose room is in the safety region of $\PP^{\mathbf{n}}$ is true.

We aim at studying how the safety bounds depend on the dimensions of the factors of the considered Segre variety. In order to reach this target, we prove the following result.

\begin{teo}\label{teo: stimaperO+}
 Let $X = \PP^{n_1} \ttimes \PP^{n_k}$ be a Segre variety, such that $n_1 \leq \dots \leq n_k$. Then
$$\OOO_+ \geq n_k \left( \sum_{i=1}^{k-1} n_i -1\right).$$
\begin{proof}
 We prove that there is a false statement whose room is $n_k \left( \sum_{i=1}^{k-1} n_i -1\right) -1$. Consider
$$\mathcal{T}=T(n_1 \vvirg n_k;1;0 \vvirg 0,a_k),$$
where $a_k = \prod_{i=1}^{k-1} (n_i +1) - \sum_{i=1}^{k-1}n_i$. Let $\RRR$ denote its room. By computation, we can easily note that $\RRR = n_k \left( \sum_{i=1}^{k-1} n_i -1\right) -1$.

Let $X$ be the Segre variety $\PP^{n_1} \ttimes \PP^{n_k}$ and consider $1+a_k$ general points on $X$, $p,q_1 \vvirg q_{a_k}$ with the following shape:
\begin{equation*}
\begin{split}
   p&= v_1 \ootimes v_k,\\
 q_j &= w_1 ^{j} \ootimes w_k ^j \quad (\text{for } j=1 \vvirg a_k).
  \end{split}
\end{equation*}

Suppose that $\mathcal{T}$ is true. Thus the space $L = T_p X + \sum_{j=1}^{a_k} (T_{q_j}X)_k$ has the expected dimension, that is:
$$ \dim L = \left(1+ \sum_{i=1}^{k-1} n_i \right) + a_k (n_k+1).$$
We have that $\dim T_p X = 1 + \sum_1 ^k n_i$, $\dim (T_{q_j}X)_k = n_k +1$ and then $\dim \left(\sum_{j=1}^{a_k} (T_{q_j}X)_k \right) \leq a_k(n_k+1)$. Thus, if $\mathcal{T}$ is true, $T_p X$ and $\left(\sum_{j=1}^{a_k} (T_{q_j}X)_k \right)$ intersect properly. We show that these spaces have two subspaces that do not intersect properly. Consider

\begin{equation*}
\begin{split}
T_{\widehat{k}} &= V_1 v_2 \dots v_{k-1} v_k + \dots + v_1 \dots v_{k-2} V_{k-1} v_k , \\
 Wv_k &= \langle w_1^j \dots w_{k-1}^j v_k \vert j=1 \vvirg a_k \rangle .
\end{split}
\end{equation*}

Since the points are in general position, we have that $\dim Wv_k = a_k$. Besides $\dim T_{\widehat{k}} = 1 + \sum_1 ^{k-1} n_i$. Since both are subspaces of $V_1 \dots V_{k-1}v_k$, by Grassmann's formula, we have
\begin{equation*}
\begin{split}
 \dim [ Wv_k \cap T_{\widehat{k}} ] &\geq \dim(W v_k)+ \dim (T_{\widehat{k}} ) - \dim (V_1 \dots V_{k-1} v_k) = \\
&= a_k +  1 + \sum_{i=1}^{k-1} n_i - \prod_{i=1} ^{k-1} (n_i +1) = 1.
\end{split}
\end{equation*}

Since these two spaces do not intersect properly, $\mathcal{T}$ is false.
\end{proof}
\end{teo}

Moreover, we can give a more precise result if we consider the product of copies of $\PP^1$:

 \begin{teo}\label{teo: O+perP1}
 Let $X= \left(\PP^1\right)^k$ with $k\geq 3$ and let $\OOO_+$ be its positive safety bound. Then $\OOO_+ \in \{ k-2,k-1\}$.
\begin{proof}
By Theorem \ref{teo: stimaperO+}, we have that $\OOO_+ \geq k-2$. We will prove by induction on $k$ that any statement $\mathcal{T}$ with room $\RRR \geq k-1$ is true. The cases $k=3,4,5,6,7$ can be proved by explicit computation (see Section \ref{section:resultsobtained}): we can check that in these cases $\OOO_+$ is exactly $k-2$. This provides the base of the induction. Moreover, since $\left(\PP^1\right)^k$ has no defective secant varieties if $k\geq 5$ (see \cite{CGG:prodottiP1}), we may assume that at least one of the $a_i$ is non-zero (w.l.o.g. assume $a_1\neq 0$).

Thus, let $k\geq 8$ and let $\mathcal{T}=T(1\vvirg 1;s;a_1 \vvirg a_k)$ be a statement of room $\RRR \geq k-1$. We may consider $\RRR \in \{k-1,k \}$, since if the room is greater we can prove that $T(1 \vvirg 1;s;b,a_2 \vvirg a_k)$, with $b \geq a_1$ such that its room is in $\{k,k-1 \}$, is true, and conclude by Remark \ref{remark: ordinestatements}. Moreover note that condition \eqref{eq:condizioneroomproporzionataammissibile} holds for $\gamma'=\gamma'' = \frac{1}{2}$, so, for any $j=1\vvirg k$, there exists an eligible reduction on the $j$-th factor (and in particular for this reduction \eqref{eq:condizioneroomproporzionata} holds).

First, suppose that $a_1 \geq k-2$. In this case, apply a reduction on the first factor, obtaining two subabundant statements $\mathcal{T}'$ and $\mathcal{T}''$, with room $\RRR'$ and $\RRR''$ respectively, concerning the Segre variety $\PP^0 \times \left(\PP^1\right)^{k-1}$. By Proposition 2.6 these statements are true if and only if the corresponding statements concerning $\left( \PP^1 \right)^{k-1}$, are true. Let $\mathcal{T}'_*,\mathcal{T}''_*$ be the statements concerning $\left( \PP^1 \right)^{k-1}$. Note that the room of $\mathcal{T}'_*$ is at least $\RRR'+a_1 \geq k-2$ and so the statement is true by the induction hypothesis. Similarly, $\mathcal{T}''_*$ is true. So, if one of the $a_i$ is greater or equal then $k-2$, $\mathcal{T}$ is true.

Consider $a_i \leq k-3$ for any $i$. Consider any eligible reduction of $\mathcal{T}$ on the first factor in two subabundant statements $\mathcal{T}',\mathcal{T}''$, with room $\RRR',\RRR''$ respectively. We prove that $s',s'' \geq k-3$. Suppose $s' \leq k-4$: since
$$\RRR' = 2^{k-1}  - s' k - a_1 - 2 \textstyle \sum_{i\neq 1} a_i,$$
we obtain
$$\RRR' \geq 2^{k-1} - (k-4)k - (k-3) - 2(k-1) (k-3)  = 2^{k-1} - [3k^2-11k+3].$$
Since the reduction is eligible, we have $\RRR' \leq \RRR \leq k$ and so
$$2^{k-1} - [3k^2-11k+3] \leq k.$$
This inequality does not hold if $k \geq 8$, so $s' \geq k-3$. Similarly $s'' \geq k-3$.

Similarly to the first part of the proof, we obtain statements concerning $\left(\PP^1\right)^{k-1}$, and, since $a_1 \geq 1$ and $s',s'' \geq k-3$, their room is greater then $k-2$. By the induction hypothesis, we have that $\mathcal{T}$ is true.
\end{proof}
 \end{teo}

\section{Control of room value and growth of safety bounds}

In this section we will provide the tools that are necessary to prove Theorem \ref{Thm: stima asintotica per la non difettivita}. First, we will use Theorem \ref{teo: room proporzionata} iteratively, in order to estimate, in a particular case, the safety bounds of a Segre variety, by supposing to know the safety bounds of some low dimensional Segre varieties. Actually, in the next section we will explicitly compute the safety region of some low dimensional Segre varieties.

First, we claim the following remark, whose proof is immediate:

\begin{remark}\label{remark: ammissibilita dalle foglie}
 Consider a statement $\mathcal{T}=T(n_1 \vvirg n_k;s;a_1 \vvirg a_k)$. Let $\mathfrak{R}$ be a reduction tree for $\mathcal{T}$ and let $\mathcal{T}_1 \vvirg \mathcal{T}_r$ be the statements corresponding to the leaves of the tree. Then, all the reductions of $\mathfrak{R}$ are eligible if and only if $\mathcal{T}_1 \vvirg \mathcal{T}_r$ have the same abundance.
\end{remark}

This remark will be fairly useful in this section because we will find reduction trees without caring for the eligibility of the reductions: at the end of the process we find conditions guaranteeing that the leaves of our tree have the same abundance and then we conclude by Remark \ref{remark: ammissibilita dalle foglie}.

Initially, we deal with a particular case, which allows us to apply \emph{central} reductions, namely reductions that generate statements referring to the same Segre variety or, in other words, in which the value of the dimension of the factor on which the statement is reduced is divided in two equal integers. These results will allow us to estimate the positive safety bound of some particular Segre variety. By some simple inequalities, this leads to the proof of Theorem \ref{Thm: stima asintotica per la non difettivita}.

Consider a Segre variety $X = \PP^{n_1} \ttimes \PP^{n_k}$, with $n_j +1 = 2^{d_j}c$, $d_1 \leq \dots \leq d_k$. Given a statement $\mathcal{T}=T(n_1 \vvirg n_k ; s;a_1 \vvirg a_k)$, of room $\RRR$, we explain our strategy splitting it in several steps.

The first step consists in reductions on the $k$-th factor, in order to lead its dimension to be $c-1$. So we have $d_k$ sub-steps: in each one of them, we apply a central reduction on the $k$-th factor of all the statements obtained after the previous sub-step (or simply applying a reduction on $\mathcal{T}$ for the first sub-step). We will obtain $2^{d_k}$ statements of the form

$$T(n_1\vvirg n_{k-1},c-1;s^*;a_1^*,\vvirg ,a_k^*).$$

The next steps consist in following the same procedure with each one of the factors, applying the reductions on the factor of greatest dimension. At the end of the whole process, we obtain $2^{d_1+\dots + d_k}$ statements, concerning the Segre variety $Y=\PP^{c-1} \ttimes \PP^{c-1}$. Supposing to know the safety region of $Y$, we aim to estimate the safety bounds of $X$, by applying Theorem \ref{teo: room proporzionata}. 

\begin{lem}\label{lemma: 2dj->c}
 Let $\mathcal{T}=T(n_1 \vvirg n_k;s;a_1 \vvirg a_k)$ be a statement with room $\RRR$, with $n_1+1=2^{d}c$. Then there exists a reduction tree for $\mathcal{T}$, whose leaves are $2^{d}$ statements concerning the Segre variety $\PP^{c-1}\times \PP^{n_2}\ttimes \PP^{n_k}$, with room generically denoted by $\RRR_{d}$ such that
\begin{equation}\label{eq: lemma 2dj->c : iniziale}
\left\vert \RRR_d- \frac{1}{2^d} \RRR \right\vert \leq \textstyle\sum_2^{k} n_i.
\end{equation}

\begin{proof}
 Denote by $\mathcal{T}_1$ the two statements obtain after a central reduction on the first factor of $\mathcal{T}$ and denote generically by $\RRR_1$ their room. By Theorem \ref{teo: room proporzionata}, we have 
\begin{equation}\label{eq: lemma 2dj->c : baseinduzione}
 \left\vert \RRR_1 - \frac{1}{2}\RRR \right\vert \leq \frac{1}{2} \left( \textstyle\sum_2^{k} n_i\right). 
\end{equation}
Denote by $\mathcal{T}_{\ell}$ whichever of the statements obtained after a central reduction on the first factor of one of the $\mathcal{T}_{\ell -1}$, and generically by $\RRR_\ell$ its room. By induction it is easy to check that
\begin{equation}\label{eq: lemma 2dj->c : intermedia}
 \left\vert \RRR_\ell - \frac{1}{2^\ell} \RRR \right\vert \leq \left(\textstyle\sum_2^{k} n_i\right) \left(\textstyle\sum_1 ^{\ell} 2^{-j}\right).
\end{equation}

Indeed, the base of the induction is given by \eqref{eq: lemma 2dj->c : baseinduzione}, while if $\ell \geq 2$ we have
\begin{equation*}
\begin{split}
  \RRR_{\ell} &\leq \frac{1}{2} \RRR_{\ell-1} + \frac{1}{2} \textstyle \sum_{2}^{k} n_i \leq \\
&\leq \frac{1}{2}\left[ \frac{1}{2^{\ell-1}} \RRR + \left(\textstyle\sum_2^{k} n_i\right) \left(\textstyle\sum_1 ^{\ell-1} 2^{-j}\right)\right]+ \frac{1}{2} \textstyle \sum_{2}^{k} n_i \leq\\
&\leq \frac{1}{2^\ell}\RRR+\left(\textstyle\sum_2^{k} n_i\right) \left(\textstyle\sum_1 ^{\ell} 2^{-j}\right),
\end{split}
\end{equation*}
where the first passage is given by Theorem \ref{teo: room proporzionata} and the second one by the induction hypothesis. Similarly we obtain
$$\RRR_{\ell}\geq  \frac{1}{2^\ell}\RRR-\left(\textstyle\sum_2^{k} n_i\right) \left(\textstyle\sum_1 ^{\ell} 2^{-j}\right),$$
and then \eqref{eq: lemma 2dj->c : intermedia}.

So, for $\ell = d$, we have
$$ \left\vert \RRR_d- \frac{1}{2^d} \RRR \right\vert \leq \left(\textstyle\sum_2^{k} n_i\right) \left(\textstyle\sum_1 ^{d} 2^{-j}\right).$$

To conclude, consider that $\sum_1 ^\infty 2^{-j} = 1$.\end{proof}
\end{lem}

Now suppose to iterate the procedure on each factor; we have the following result:
\begin{prop}\label{prop: n1..nk -> c..c}
 Let $\mathcal{T}=T(n_1 \vvirg n_k;s;a_1 \vvirg a_k)$ be a statement with room $\RRR$, such that $n_j+1=2^{d_j}c$ for any $j=1 \vvirg k$, with $d_1 \geq \dots \geq d_k$. Then there exists a reduction tree for $\mathcal{T}$, whose leaves are $2^{d_1+\dots + d_k}$ statements concerning the Segre variety $Y=\left(\PP^{c-1}\right)^k$, with room generically denoted by $\RRR_{fin}$ such that
\begin{equation}\label{eq: prop: n1..nk -> c..c - tesi}
 \left\vert \RRR_{fin} - \frac{1}{2^{d_1 + \dots + d_k}} \RRR\right\vert \leq  c(k-1)^2 - \frac{k(k-1)}{2}.
\end{equation}

\begin{proof}
 By Lemma \ref{lemma: 2dj->c}, we can reduce $\mathcal{T}$ in $2^{d_1}$ statements $\mathcal{T}_1$, with room generically denoted by $\RRR_1$, such that \eqref{eq: lemma 2dj->c : iniziale} holds. Each one of the statements $\mathcal{T}_1$ can be reduced, on the second factor, obtaining statements concerning $\PP^{c-1}\times \PP^{c-1} \times \PP^{n_3}\ttimes \PP^{n_k}$, whose room is generically denoted by $\RRR_2$ and such that:

\begin{equation*}
 \left\vert \RRR_2 - \frac{1}{2^{d_2}} \RRR_1 \right\vert \leq (c-1) + \left(\textstyle\sum_{3}^{k} n_i \right).
\end{equation*}

More generally, for any $\ell=1 \vvirg k$ denote by $\mathcal{T}_{\ell}$ any statement obtained after $d_\ell$ reductions on the $\ell$-th factor of $\mathcal{T}_{\ell-1}$. We have
$$ \left\vert \RRR_{\ell} - \frac{1}{2^{d_\ell}} \RRR_{\ell-1}\right\vert \leq (\ell-1)(c-1) + \left( \textstyle \sum_{\ell+1}^k n_i\right).$$

By induction we prove

\begin{equation}\label{eq:prop: n1..nk -> c..c - intermedia}
\left\vert \RRR_\ell - \frac{1}{2^{d_1 + \dots + d_\ell}} \RRR \right\vert \leq \sum_{t=2} ^ {\ell+1} \frac{1}{2^{d_t +..+d_\ell}} \left(\textstyle\sum_t^k n_i \right) + \sum_{t=2}^{\ell} \frac{1}{2^{d_{t+1} + \dots + d_\ell}} (t-1)(c-1).
\end{equation}

The base of the induction is given by Lemma \ref{lemma: 2dj->c}. If $\ell \geq 2$, we have
\begin{equation*}
\begin{split}
\RRR_\ell &\leq \frac{1}{2^{d_\ell}} \RRR_{\ell -1}+ (\ell -1)(c-1) + \left( \textstyle\sum _{\ell+1}^k n_i\right) \leq \\
&\leq \frac{1}{2^{d_\ell}} \left[ \frac{1}{2^{d_1 + \dots + d_{\ell-1}}} \RRR + \sum_{t=2} ^ {\ell} \frac{1}{2^{d_t +..+d_{\ell-1}}} \left(\textstyle\sum_t^k n_i \right) + \sum_{t=2}^{\ell-1} \frac{1}{2^{d_{t+1} + \dots + d_{\ell-1}}} (t-1)(c-1)\right] +\\
& \hspace{.6\textwidth} +  (\ell -1)(c-1) + \left( \textstyle\sum _{\ell+1}^k n_i\right) \leq \\
& \leq  \frac{1}{2^{d_1 + \dots + d_\ell}} \RRR + \sum_{t=2} ^ {\ell+1} \frac{1}{2^{d_t +..+d_\ell}} \left(\textstyle\sum_t^k n_i \right) + \sum_{t=2}^{\ell} \frac{1}{2^{d_{t+1} + \dots + d_\ell}} (t-1)(c-1).
\end{split}
\end{equation*}
 Similarly, we also have
$$\RRR_\ell \geq \frac{1}{2^{d_1 + \dots + d_\ell}} \RRR + \sum_{t=2} ^ {\ell+1} \frac{1}{2^{d_t +..+d_\ell}} \left(\textstyle\sum_t^k n_i \right) + \sum_{t=2}^{\ell} \frac{1}{2^{d_{t+1} + \dots + d_\ell}} (t-1)(c-1),$$
and so \eqref{eq:prop: n1..nk -> c..c - intermedia}.
If $\ell = k$, since $n_j+1= 2^{d_j}c$, the following estimate holds:

\begin{equation*}
\begin{split}
\sum_{t=2}^{k+1} \frac{1}{2^{d_t + \dots + d_k}}\left( \textstyle \sum_t^k n_i \right) &= \sum_{t=2}^{k} \left({\textstyle \sum}_{t}^k \left(\frac{2^{d_i}c -1}{2^{d_{t} + \dots +d_{k}}} \right)\right) \leq  \\
&\leq \sum_{t=2}^k \left({\textstyle\sum}_t ^k c \right) \leq \sum_{t=2}^k c(k-t) = c \frac{(k-1)(k-2)}{2}.
\end{split}
\end{equation*}
Moreover
$$\sum_{t=2}^k \frac{1}{2^{d_1 + \dots +d_k}}(t-1)(c-1) \leq (c-1)\sum_1 ^{k-1} t = (c-1)\frac{k(k-1)}{2}.$$

So we obtain
$$\left\vert \RRR_{fin} - \frac{1}{2^{d_1 + \dots + d_k}} \RRR\right\vert \leq  c \frac{(k-1)(k-2)}{2}+(c-1)\frac{k(k-1)}{2} = c(k-1)^2 - \frac{k(k-1)}{2}.$$

\end{proof}
\end{prop}

\begin{teo}\label{teo: sicurezza d1..dk}
 Let $X = \PP^{n_1} \ttimes \PP^{n_k}$ be a Segre variety, with $n_j +1 = 2^{d_j} c$ and $k\geq 3$ and let $Y = \left(\PP^{c-1}\right)^k$. Let $\OOO_+ , \OOO_-$ be the safety bounds of $Y$ and $\OOO_+ ^* , \OOO_- ^*$ be the safety bounds of $X$. Then 
\begin{equation}
 \begin{split}
  \OOO_+ ^* &\leq 2^{d_1 + \dots +d_k}(\OOO_+ + C), \\
  \OOO_- ^* &\geq 2^{d_1 + \dots +d_k}(\OOO_- - C),
 \end{split}
\end{equation}
where  $C=c(k-1)^2 - \frac{k(k-1)}{2}$. %$C = \left( 2^{-2d_2} + 2^{-d_1-d_2} + 2^{-2 d_1 +1} + 2^{d_2 - 3d_1 -1} +3 \right) c$.
\begin{proof}
Consider a statement $\mathcal{T}=T(n_1\vvirg n_k;s;a_1\vvirg a_k)$, of room $\RRR \geq  2^{d_1+d_2+d_3}(\OOO_+ + C)$. We show that $\mathcal{T}$ is true. Proposition \ref{prop: n1..nk -> c..c} provides a reduction tree for $\mathcal{T}$ whose leaves are statements $\mathcal{T}_{fin}$ concerning the variety $Y$, with room generically denoted by $\RRR_{fin}$ such that \eqref{eq: prop: n1..nk -> c..c - tesi} holds. Hence
$$\RRR_{fin} \geq \frac{1}{2^{d_1 +\dots +d_k}} \RRR -C \geq \frac{1}{2^{d_1 +\dots + d_k}} \left[ 2^{d_1+\dots +d_k}(\OOO_+ + C) \right] -C \geq \OOO_+.$$
Thus, all the leaves are subabundant true statements. Thus, by Remark \ref{remark: ammissibilita dalle foglie}, we have that all the reductions are eligible, and hence, $\mathcal{T}$ is true.
In similar manner, by considering a superabundant statement, we have the other inequality.
\end{proof}
\end{teo}

Theorem \ref{teo: sicurezza d1..dk} and Theorem \ref{teo: O+perP1} provide the proof of Theorem \ref{Thm: stima asintotica per la non difettivita}.

\section{Proof of Theorem \ref{Thm: stima asintotica per la non difettivita}}

\begin{PROOF1.1}
First, consider a Segre variety of the shape $X=\PP^{n_1} \ttimes \PP^{n_k}$, with $n_j+1 = 2^{d_j+1}$ and let $\OOO_+$ be its positive safety bound. Let $\mathcal{T}=T(n_1 \vvirg n_k;s;0\vvirg 0)$ be the statement describing the $s$-secant variety of $X$ and let $\RRR$ be its room. By Theorem \ref{teo: O+perP1}, the positive safety bound of $\left(\PP^1 \right)^k$ is lower than $k-1$, so by Theorem \ref{teo: sicurezza d1..dk} (with $c=2$), we have
$$\OOO_+ \leq 2^{d_1 + \dots + d_k}(k-1 + 2(k-1)^2 - \frac{k(k-1)}{2}).$$
If $\RRR$ is greater than the value on the right side of the inequality, we can deduce that the $s$-secant variety of $X$ is non-defective. So, since $2^{d_1 + \dots +d_k} = \frac{\prod (n_i+1)}{2^k}$, we have that $T(n_1\vvirg n_k;s;0^k)$ is true if

$$ \textstyle \prod (n_i+1) - s (1+\sum n_i) \geq \dfrac{1}{2^k}\prod(n_i+1)\left(\dfrac{3k^2-5k+2}{2}\right).$$

Easily we obtain that an integer $s$ satisfies this inequality if and only if it satisfies the following one:

$$ \textstyle s \leq \dfrac{\prod(n_i+1)}{1+\sum n_i} \left( 1 - \dfrac{3k^2-5k+2}{2^{k+1}} \right).$$

Let $\Delta_k =  1 - \dfrac{3k^2-5k+2}{2^{k+1}}$. Note that for any $k$ the value $\Delta_k$ is positive.

Now, consider any Segre variety $\PP^{n_1} \ttimes \PP^{n_k}$. For any $j$, let $m_j +1 =2^{\lfloor \log _2 (n_i+1) \rfloor}$ and let $Y = \PP^{m_1} \ttimes \PP^{m_k}$. $Y$ is of the shape dealt in the previous part of the proof, and so its $s$-secant variety is non defective if
$$ s \leq \frac{\prod(m_i+1)}{1+\sum m_i} \Delta_k.$$

Since $m_j \leq n_j$ for any $j$, by Remark \ref{remark: ordinestatements}, for these values of $s$, the $s$-secant varieties of $X$ are also non-defective. So, we have that $T(n_1 \vvirg n_k;s;0^k)$ is true if
$$ s \leq \frac{\prod(n_i+1)}{1+\sum n_i} \ \frac{\prod(m_i+1)}{\prod(n_i+1)} \ \frac{1+\sum n_i}{1+\sum m_i} \Delta_k.$$

Since $\frac{1}{2} (n_i+1) \leq m_i+1 \leq (n_i+1)$, we have that this condition holds if
$$s \leq \frac{\prod(n_i+1)}{1+\sum n_i} \frac{\Delta _k}{2^k}.$$
Set $\Theta_k = \frac{\Delta _k}{2^k}$ and conclude. $\hfill \qed$
 \end{PROOF1.1}

The value $\Theta_k$ provides a bound for the ratio between the greater $s$ such that the $s$-secant variety of a Segre variety with $k$-factors is non-defective and it expected rank. However, for secant varieties of a particular shape, this bound can be strongly improved:

\begin{corol}\label{cor: Teo 1.1}
 Let $X = \PP^{n_1} \ttimes \PP^{n_k}$ be a Segre variety with $n_i +1 =2^{d_i+1}$. Then, the $s$-secant variety of $X$ is non-defective for any $s$ such that
$$ s \leq \Delta_k \frac{\prod(n_i+1)}{1+\sum n_i},$$
where $\Delta_k = 1 - \dfrac{3k^2-5k+2}{2^{k+1}}$.
\end{corol}

\section{Results obtained by explicit computation}\label{section:resultsobtained}

An algorithm implemented in Macaulay2 allows us to compute by a Monte-Carlo method whether a given statement $T(\mathbf{n},s,\mathbf{a})$ is true. The script considers $s + \sum a_i$ random points of $X=\PP^{n_1}\ttimes \PP^{n_k}$ and computes the dimension of the space $L$ of Definition \ref{def:Tnsa}. If this dimension equals the expected dimension, by a semicontinuity argument we can deduce that the statement is true. If the dimension of the space is lower than the expected one, we cannot conclude that the statement is false, since the random choice made by the script could have led to a singular situation. However, in these cases, the answer of the script is a really strong clue that the statement is false.

Using this script in a loop makes it possible to compute the safety region of low-dimensional Segre variety. We computed the safety region of Segre varieties with three factors whose dimension is lower than $5$ and of Segre varieties with four factors whose dimension is lower than $3$. We saw that the positive and negative safety bounds in these cases are one opposite to the other, and that the positive one is equal to the lower-bound that we gave in Theorem \ref{teo: stimaperO+}. Moreover we computed the positive safety bound of $\left(\PP^1\right)^5 , \left(\PP^1\right)^6$ and $\left(\PP^1\right)^7$, which are $3,4$ and $5$ respectively.

In the following table, we sum up all the results we have about the safety bounds:

\begin{center}
\begin{tabular}{ccc}
$\PP^{n_1} \ttimes \PP^{n_k}$ & $\OOO_+ \geq \left(\sum_{i\leq k} n_i -1\right)n_k$ & Thm. \ref{teo: stimaperO+}
\\
\hline
$\left( \PP^1 \right) ^k$ & $k-2 \leq \OOO_+ \leq k-1$ & Thm. \ref{teo: O+perP1}\\
\hline
\begin{tabular}{c}$\PP^{n_1}\times \PP^{n_2} \times \PP^{n_3}$\\ with $n_1 \leq n_2 \leq n_3 \leq5$\end{tabular} & \begin{tabular}{c}$\OOO_+ =(n_1 + n_2 -1)n_3$ \\$\OOO_- = -(n_1 +n_2 -1)n_3$\end{tabular} & by the script\\
\hline
\begin{tabular}{c}$\PP^{n_1}\times \PP^{n_2} \times \PP^{n_3} \times \PP^{n_4}$\\ with $n_1 \leq n_2 \leq n_3 \leq n_4 \leq3$\end{tabular} & \begin{tabular}{c}$\OOO_+ =(n_1 + n_2 + n_3 -1)n_4$ \\$\OOO_- = -(n_1 +n_2 + n_3-1)n_4$\end{tabular}& by the script\\
\hline
$(\PP^1)^k$ with $k = 5,6,7$ & $\OOO_+ = k-2$& by the script
\end{tabular}
\end{center}

This suggests the following conjecture:

\begin{congettura}
 Let $X = \PP^{n_1} \ttimes \PP^{n_k}$ be a Segre variety with $n_1 \leq \dots \leq n_k$ and let $\OOO_-$ and $\OOO_+$ be the safety bounds of $X$. Then the following statements hold:

\begin{itemize}
\item[-] there is a symmetry between the two safety bounds, i.e. $\OOO_- = - \OOO_+$;
\item[-] the positive safety bound is $\OOO_+ = \left(\sum_{i=1}^{k-1} n_i -1 \right)n_k$;
\item[-] if $a_k = \prod_{i=1}^{k-1} (n_i +1) - (n_k+1) \geq 0$, the statement $T(n_1 \vvirg n_k;n_k; 0 \vvirg 0, a_k)$, of room $\RRR$, satisfies $\RRR = \OOO_-$ and is false.
\end{itemize} 
\end{congettura}

Moreover, we analyzed the secant varieties of Segre varieties with four factors of the same dimension, namely of the shape $\PP^n \times \PP^n\times \PP^n\times \PP^n$, for $n\leq 10$, finding some interesting results. With the aid of the script, we can prove the following theorem:

\begin{teo}
 Let $X$ be a Segre variety of the shape $\left( \PP^n \right)^4$, with $2 \leq n \leq 10$. Then, $X$ has no defective secant varieties except at most in the two following cases
\begin{itemize}
\item[-] the $199$-secant variety of $\left(\PP^8\right)^4$;
\item[-] the $357$-secant variety of $\left(\PP^{10}\right)^4$.
\end{itemize}
\end{teo}
The proof of this result is provided by the script for $n\leq 6$. As for secant varieties of $\left(\PP^n\right)^4$, with $n=7,8,9,10$, by computing their reduction trees (see \cite{TesiMia}, Section 6.1), we proved that they are all nondefective, except at most in the two listed cases, that we cannot approach by the inductive technique, because of arithmetical reasons. As far as $n=1$ is concerned, it is well known that the $3$-secant variety of $\left(\PP^1\right)^4$ is defective (see \cite{aboottapet} or \cite{CGG:prodottiP1}). 

Actually, even though the statements corresponding to the two listed secant varieties can not be managed by the inductive technique, it is unlikely that they are defective; indeed, even the $96$-secant variety of $(\PP^6)^4$ cannot be managed by the inductive technique, but in this case the script shows that it is nondefective.

\end{document}